\documentstyle{amsart}

\newtheorem{theorem}{Theorem}[section]
\newtheorem{defn}[theorem]{Definition}
\newtheorem{prop}[theorem]{Proposition}
\newtheorem{lem}[theorem]{Lemma}
\newtheorem{cor}[theorem]{Corollary}

\newcommand{\sub}{\subseteq}

\newcommand{\sm}{\setminus}
\newcommand{\proof}{\noindent {\bf Proof: }}
\newcommand{\noi}{\noindent}
\newcommand{\dom}{\mbox{ dom}}

\newsymbol\upharpoonright 1316
\newcommand{\rest}{\upharpoonright}

\newcommand{\ra}{\rangle}
\newcommand{\la}{\langle}

\newcommand{\CU}{{\cal U}}
\newcommand{\CV}{{\cal V}}
\newcommand{\CX}{{\cal X}}

\newcommand{\al}{\alpha}
\newcommand{\be}{\beta}
\newcommand{\de}{\delta}
\newcommand{\ga}{\gamma}
\newcommand{\ifiner}{${\frak if}$\/}
\newcommand{\eifiner}{${\frak eif}$\/}

\newcommand{\fo}{\; ^\omega  \! \omega}
\newcommand{\oo}{\; \omega^\omega}
\newcommand{\fw}{\; ^\omega  \! \omega_1}

\newcommand{\so}{[\omega]^{\omega}}

\newcommand{\ff}{\; ^{\omega} \! 2}

\newcommand{\po}{\mbox{{\Large $\wp$}}(\omega)}
\newcommand{\fa}{\forall}
\newcommand{\fai}{\forall^\infty}
\newcommand{\ex}{\exists}

\begin{document}

\title{A few special ordinal ultrafilters}

\author{Claude Laflamme}
\address{Department of Mathematics and Statistics \\
         University of Calgary \\
         Calgary, Alberta  \\
         Canada T2N 1N4}

\email{laflamme@@acs.ucalgary.ca}
\date{}
\thanks{This research was partially supported by  NSERC of Canada.}
\maketitle

\begin{abstract}
We prove various results on the notion of  ordinal ultrafiters
introduced by J. Baumgartner. In particular, we show that this notion
of ultrafilter complexity is independent of the more familiar
Rudin-Keisler ordering.
\end{abstract}
			     
\section{Introduction}

$ x= y= \mbox{ \boldmath{$ z $} } $

Interesting ultrafilters are those comprising rich combinatorial
properties of some sort. Traditional criterions consist of partition
relations on the natural numbers and the Rudin-Keisler ordering. In
\cite{Bau}, Baumgartner introduces several new combinatorial notions
for ultrafilters and we show in this paper that his concept of ordinal
ultrafilter, related to the behaviour of functions from $\omega$ to
$\omega_1$, is independent of the traditional combinatorics and
therefore brings a new insight in the theory of ultrafilters. 

Our terminology is standard but we review the main concepts and
notation.  The natural numbers will be denoted by $\omega$, $\ff$ and
$\fo$ denote the collection of functions from $\omega$ to 2 and to
$\omega$ respectively; similarly, $\po$ and $\so$ denote the
collection of all and infinite subsets respectively.  We can view
members of $\po$ as members of $\ff$ by considering their
characteristic functions.  

A filter is a collection of subsets of $\omega$ closed under finite
intersections, supersets and to avoid trivialitites contain all
cofinite sets; it is called proper if it contains only infinite sets.
Given a collection $\CX \sub \po$, we let $\la \CX \ra$ denote the
filter generated by $\CX$. An ultrafilter is a proper maximal filter.

Here are a few examples of combinatorially rich ultrafilters (see
\cite{Boo}).

\begin{defn}
An ultrafilter $\CU$ is called a

\begin{enumerate}
\item  P-point if for any $f \in \fo$, there is an $X \in
\CU$ such that $f \rest X$ is either constant or finite-to-one.
\item Ramsey ultrafilter if $\CU$ contains a homogeneous set  for each
$f:[\omega]^k \rightarrow \ell$, $k,\ell \in \omega$. 
\end{enumerate}
\end{defn}

The well-known Rudin-Keisler ordering for ultrafilters is defined by
\[ \CU <_{RK} \CV \mbox{ if } (\ex f \in \fo) \CU = \la \{ f''X : X \in
\CV \} \ra. \]

There are some important connections between the previous notions,
indeed $\CU$ is a Ramsey ultrafilter if and only if it is minimal in
the Rudin-Keisler ordering, see \cite{Boo} for more. 

We recall the basic operations of multiplication and exponentiation on
ordinals.

\begin{defn}For any ordinals  $\alpha, \beta$,
\begin{enumerate}
\item $\alpha \cdot 0 = 0$
\item $\alpha \cdot 1 = \alpha$
\item $\alpha \cdot (\beta + 1) = \alpha \cdot \beta + \alpha$
\item If $\beta$ is a limit ordinal,
     then $\alpha \cdot \beta = \sup \{ \alpha \cdot \xi: \xi < \beta \}.$
\item $\alpha^0 = 1$
\item $\alpha ^{\beta + 1} = \alpha^\beta + \alpha$
\item If $\beta$ is a limit ordinal,
      then $\alpha^\beta = \sup \{ \alpha^\xi : \xi < \beta \}$.
\end{enumerate}
\end{defn}

\noi As any subset $X$ of ordinals is well ordered, we can define the order
type of $X$ as the unique ordinal order isomorphic to $X$.

\section{Basic ordinal ultrafilters}

We recall Baumgartner's notion of ordinal ultrafilter and a few
related tools. 

\begin{defn} Let $\alpha \leq \omega_1$ be any ordinal and $\CU$ an
ultrafilter on $\omega$.
 \begin{enumerate}
\item $\CU$ is said to be an $\alpha$-ultrafilter if $\alpha$ is the
smallest ordinal such that for every $h:\omega \rightarrow \omega_1$
we can find an $X \in \CU$ such that $h '' X$ has order type at most
$\alpha$.
\item $\CU$ is a strict $\alpha$-ultrafilter if  in the above
definition we demand that the order type of $h '' X$ is strictly less
than $\alpha$. 
\item (The infinite Rudin-Keisler ordering) $\CU <_\infty \CV$ if
there is $f \in \fo$ with $f(\CV)=\CU$ (so $\CU <_{RK} \CV$) but $f
\rest X$ is not finite-to-one or constant for any $X \in \CV$.
\end{enumerate}
\end{defn}

Here are some basic known results on ordinal ultrafilters.

\begin{prop}
(Baumgartner \cite{Bau})
\begin{enumerate}
\item If $\CU$ is an $\alpha$-ultrafilter, then $\alpha$ is an
indecomposable ordinal, that is $\alpha = \omega^\beta$ for some
$\beta$. 
\item P-points are exactly the $\omega$-ultrafilters.
\end{enumerate}
\end{prop}

Indeed, if $\CU$ is an ultrafilter and $h \in \fw$, then
\[ \min \{ \alpha: (\ex X \in \CU) \mbox{ $h''X$ has order type
$\alpha$ } \} \]
\noi must be an indecomposable ordinal. As a generalisation of the
second result we have the following.

\begin{prop}\label{omega-k} Let $k \in \omega$ and $\CU$ an
ultrafilter such that 

\[ \mbox{(*) }(\fa h \in \fw)(\ex X \in \CU)
 \mbox{ the order type of $h''X$ is strictly less than $\oo$}. \]

\noindent  Then $\CU$ is an $\omega^{k}$-ultrafilter precisely when it has a
$<_\infty$-chain of length $k$ below it (possibly including $\CU$) but
no $<_\infty$-chain of length $k+1$.
\end{prop}

We break the proof into a few lemmas that will remain useful later for
other purposes.

\begin{lem}\label{chain}
 Let $k \in \omega$ and suppose that we have $\CU_0 >_\infty \CU_1
>_\infty \cdots >_\infty \CU_k$, an $>_\infty$-chain of length
$k+1$. Then there is a map $h:\omega \rightarrow \omega_1$ such that
the order type of $h'' X$ is at least $\omega^{k+1}$ for any $X \in
\CU_0$.
\end{lem}

\proof We prove the result by induction on $k$; the case $k=0$ being
obvious. 

\noi Assuming the result for $k$, consider a chain of the form 
 $\CU_0 >_\infty \CU_1 >_\infty \cdots >_\infty \CU_{k+1}$. By
induction, there is a map $g$ such that the order type of $g '' X$
is at least $\omega^{k+1}$ for each $X \in \CU_1$. Now fix a map $f
\in \fo$ witnessing $\CU_0 >_\infty \CU_1$, and define $h \in \fw$ by

\[ h(m) = \la g(f(m)), m \ra \sub \omega_1 \times \omega, \]

\noi where $\omega_1 \times \omega$ is equipped with the lexicographic
ordering. For $X \in \CU_0$, we may assume that $f^{-1}\{n\} \cap X$
is infinite for all $n \in f''X$ and since the order type of $g''f''X$
is at least $\omega^{k+1}$ by assumption, the order type of $h''X$ is
at least $\omega \cdot \omega^{k+1} = \omega^{k+2}$. The required map
with range in $\omega_1$ can now  easily be obtained. \qed

And for the other direction we have.

\begin{lem}\label{shrink}
Let $\CU$ be an ultrafilter and $h \in \fw$. If
\[k = \min \{\alpha : (\ex X \in \CU) \mbox{ $h''X$ has order type at
most $\omega^\alpha$, } \} \in \omega, \]
\noi then there is an  $<_\infty$-chain below (and including) $\CU$ of length
 $k$.
\end{lem}

\proof Fix such an ultrafilter $\CU$, a map $h \in \fw$ and $k \in
\omega$ as above. Choose $X \in \CU$ such that the order type of
$h''X$ is $\omega^k$. Let $ot:h''X \rightarrow \omega^\ell$ be the
unique order preserving bijection and we may now work with the
ultrafilter $\CV=ot(h(\CU))$ and to simplify notation we work with
ultrafilters on $\omega^k$.

\noi For $i<k-1$ we define functions $g_i:\omega^k \rightarrow
\omega^k$ by

\noi $\begin{array}{ll}
g_i(\alpha) = & \omega^{k-1}\cdot m_1 + \cdots +
                                          \omega^{i+1}\cdot m_{k-i-1} \\
              & \mbox{ if }
     \end{array}$

$
 \omega^{k-1}\cdot m_1 + \cdots + \omega^{i+1}\cdot m_{k-i-1}
 \leq \alpha <
 \omega^{k-1}\cdot m_1 + \cdots + \omega^{i+1}\cdot (m_{k-i-1}+1).
$

\noi Then we obtain
\[ \CV_0 = \CV >_{RK} 
              \CV_1 = g_0(\CV_0) >_{RK} 
              \CV_2 = g_1(\CV_1) >_{RK} \cdots >_{RK}
              \CV_{k-1}= g_{k-2}(\CV_{k-2}). \] 

\noi  Now if any of the functions $g_i$ is finite-to-one when restricted to
some member $X_i$ of $\CU_i$, then the oder type of $h ''
g_i^{-1}\{X_i\}$ would be at most $\omega^{k-1 }$, a
contradiction. Thus we have obtained an $<_\infty$-chain of length $k$
below $\CU$ and the proof is complete. \qed

\bigskip

Thus by Baumgartner's result, the classical notion of P-points can be
rephrased in terms of ordinal ultrafilters, and assuming (*), the more
general notion of $\omega^k$ ultrafilter for $k \in \omega$ can be
rephrased in terms of the RK ordering. We shall see in the next
section that the assumption (*) is necessary to make this correlation,
and that actually the notion of ordinal ultrafilter is quite
independent of the RK ordering.

Assuming the Continuum Hypothesis, or more generally Martin's axiom,
it is relatively easy to construct $\omega^k$-ultrafilters for any $k
\in \omega$ (see \cite{Laf} for a general framework). In the next
section, we consider the more interesting case of $\oo$-ultrafilters.

\section{$\oo$-Ulrafilters}

We now consider the case of $\oo$-ultrafilters, where more interesting
structure occurs. We had hoped that the length of $<_\infty$-chains
below an ultrafilter as in Proposition \ref{omega-k} was a good
indication of its ordinal complexity; indeed as a Corollary to Lemma
\ref{shrink} we have:

\begin{prop}
If $\CU$ is a strict $\oo$-ultrafilter, then $\CU$ has arbitrarily
long finite $<_\infty$-chains below it.
\end{prop}

Further, similarly to Lemma \ref{chain}, a strict $\oo$-ultrafilter
cannot have an infinite descending chain.

\begin{lem}\label{infchain}
If an ultrafilter $\CU$ has an infinite decreasing
$<_\infty$-sequence below, then there is a map $f \in \fw$ such that
the order type of $f''X$ is at least $\oo$ for any $X \in \CU$.
\end{lem}

\proof Consider an infinite descending
$<_\infty$-sequence $\CU_0 >_\infty \CU_1 >_\infty \cdots$. Fix
functions $f_i \in \fo$ witnessing $\CU_i >_\infty \CU_{i+1}$. We may
assume that $f_i^{-1}\{n\}$ is infinite for each $i$ and $n \in
\omega$. We define a map $h:\omega \rightarrow \oo$ by $h = \cup_n
h_n$ as follows.  Having defined $h_0,h_1,\cdots,h_{n-1}$, choose $k_n
\notin \cup_{i<n} \dom(h_i)$, and let 
\[ \dom(h_n)= f_0^{-1}f_1^{-1}\cdots f_n^{-1} 
       \{f_n(f_{n-1}(\cdots (f_1(f_0(k_n))))) \}
       \setminus \cup_{i<n} \dom(h_i). \]

\noi Now $h_n$ is defined as any one-to-one function which respects
the following ordering on $dom(h_n)$; for $a,b \in dom(h_n)$,

\[ a \prec b \]

\noi iff for 
\[i = \min \{j: f_j(f_{j-1}(\cdots f_0(a))) =
                            f_j(f_{j-1}(\cdots f_0(b))) \},\] 

\noi we have $f_{i-1}(\cdots f_0(a)) < f_{i-1}(\cdots f_0(b))$. This
 ordering has order type exactly $\omega^{n+1}$.

\noi Now to verify that $h$ is as required, fix $X \in \CU$ and $n \in
\omega$; we show that the order type of $h '' X$ is at least
$\omega^n$. Let $X=X_0$ and more generally for $1 \leq i \leq n$ let
$X_i = f_{i-1}(\cdots (f_0(X)))$. We may assume that for each $i \leq
n$ 
\[ (\fa x \in X_i) f_i^{-1} \{f_i(x)\} \cap X_i \mbox{ is infinite.}
\]

\noi Finally if $k_m$ is such that $m \geq n$ and 

\[f_{n-2}(\cdots (f_1(f_0(k_m)))) = f_{n-2}(\cdots (f_1(f_0(x))))\]

\noi for some $x \in X$, then the order type of $h \rest (X \cap
dom(h_m))$ is exactly $\omega^{m+1}$. \qed

\vspace{.5in}

\noi {\bf Open Problem 1:} What about the corresponding influence of
{\em increasing} $<_\infty$-chains below $\CU$?

\noi Given such an ultrafilter $\CU$ with an increasing infinite
$<_\infty$-sequence 

\[ \CU >_{RK} \cdots \CU_2 >_\infty \CU_1 >_\infty \CU_0 \] 

\noi below, fix maps $g_i$ and $f_i$ witnessing $\CU >_{RK} \CU_i$ and
$\CU_{i+1} >_\infty \CU_i$ respectively. The problem is really about
the possible connections between $g_i$ and $f_i \circ g_{i+1}$, even
relative to members of $\CU$.

\noi {\bf Open Problem 2:} Can we have an ultrafilter with
arbitrarily long finite $<_\infty$-chains below $\CU$ without infinite
such chains?

\noi This looks like the most promising way to build a strict
$\omega^\omega$-ultrafilter. 

\bigskip

We now show that ordinal complexity $\omega^\omega$ is independent of
the $<_\infty$ and even the RK ordering. Theorem \ref{complex} answers
one of baumgartner's problem in \cite{Bau}.

\begin{theorem} \label{simple} (Assume CH for example, or MA, ...)
There is an $\oo$-ultrafilter whose only RK-predecessor is a Ramsey
ultrafilter. 
\end{theorem}

\begin{theorem} \label{complex}
(Assume CH for example, or MA, ...) There is an $\oo$-ultrafilter
all of whose RK-predecessors are also $\omega^\omega$-ultrafilters.
\end{theorem}
			     
The techniques used are very similar to those of \cite{Laf}; that is
we define a countably closed partial order and prove that there is
such an ultrafilter in the forcing extension. This approach somewhat
simplifies the notation but the reader will quickly realize that all
details can be carried out assuming the Continuum Hypothesis or even
Martin's Axiom. Under this last hypothesis for example, Theorem
\ref{complex} produces a descending $<_\infty$-chain of
$\oo$-ultrafilters of order type $2^{\aleph_0}$.

\begin{defn} 
\begin{enumerate}

\item An equivalence relation $E$ is said to be {\bf infinitely finer}
than $F$, written $E <_\infty F$, if each $F$ equivalence class is an
infinite union of $E$ classes. We conversely call $F$ infinitely
coarser than $E$.

\item A sequence of equivalence classes $\la E_1,E_2,\dots, E_n \ra$
is said to be {\bf infinitely finer}, or simply \ifiner, if each
$E_i <_\infty E_{i+1}$. It is said to be {\bf eventually infinitely
finer}, or \eifiner, if for all but finitely many $E_n$ equivalence
classses $C$, the sequence $\la E_1 \rest C,E_2 \rest C,\dots, E_n
\rest C \ra$ is \ifiner.

\end{enumerate}

\end{defn}

\noi Note the special role played by the last equivalence relation in
definition (2).  Observe also the following easy fact which will be
used repeatedly in the constructions.  Given an \ifiner sequence of
equivalence relations $\la E_1,E_2,\dots, E_n \ra$ on a set $X \sub
\omega$, and given a function $f \in \fo$, then we can find $Y \sub X$
such that $\la E_1 \rest Y,E_2 \rest Y,\dots, E_n \rest Y \ra$ is
still \ifiner,
 and $f \rest Y$ is either one-one, constant or else there is an $i \leq
n$ such that $f$ is constant on the $E_i \rest Y$ classes but assumes
distinct values on distinct classes.  Similarly, if $h$ is a function
from $\omega$ to $\omega_1$, then we can ensure that the order type of
$h \rest Y$ is at most $\omega^n$ (ordinal exponentiation).

\vspace{.5in}

\noi {\bf Proof of Theorem \ref{simple}} We are ready to define
our partial order.

\begin{defn}
${\Bbb P} = \{ \la \la E^i_j:j<n_i ;X_i \ra: i \in \omega \ra : E^i_0
<_\infty \cdots <_\infty E^i_{n_i-1}$ are equivalence relations on the
disjoint infinite sets $X_i \sub \omega$, and $\limsup_{i \rightarrow
\infty} n_i = \infty \}$. 
\newline For notational simplicity, we also assume
that $E^i_0$ is the finest equivalence relation, the identity, and
that $E^i_{n_i-1}$ is the coarsest equivalence relation, with only one
equivalence class.
\end{defn}

\noi We define the ordering as follows:

\[ \la \la E^i_j:j<n_i ;X_i \ra: i \in \omega \ra  \leq
       \la \la F^i_j:j<m_i ;Y_i \ra: i \in \omega \ra \]

\noi if and only if

\noi $(\fai i)(\ex k)\; \mbox{\Large [} X_i \sub Y_k$ and $(\ex \pi:n_i
\rightarrow m_k)$ increasing maps such that 

\[ E^i_j = F^k_{\pi(j)} \rest X_i \mbox{\Large ]}. \]

\begin{lem}\label{c-closed}
${\Bbb P}$ is countably closed.
\end{lem}

The proof is straightforward. More to the point we have:

\begin{lem}\label{fw}
Given $f \in \fo$, and $\la \la F^i_j:j<m_i ;Y_i \ra: i \in \omega \ra
\in {\Bbb P}$, then there is 
\[ \la \la E^i_j:j<n_i ;X_i \ra: i \in \omega \ra  \leq
       \la \la F^i_j:j<m_i ;Y_i \ra: i \in \omega \ra \]
\noi such that either:
\[ \begin{array}{ll}
     & f\rest \cup_i X_i \mbox{ is constant, } \\
 \mbox{ or } & f\rest \cup_i X_i \mbox{ is one-one, } \\
 \mbox{ or else } & f\rest X_i \mbox{ is constant for each $i$, but
takes distinct values for different $i$'s.}
  \end{array} \]

\end{lem}

\proof Fix $f \in \fo$ and $\la \la F^i_j:j<m_i ;Y_i \ra: i \in \omega \ra
\in {\Bbb P}$. We can assume, following the comments above, that for
each $i$ we have $k_i < m_i$ such that $f$ is constant on the
$F^i_{k_i}$ classes but assumes distinct values on different classes. 

\noi If $\limsup_i k_i = \infty$, then for each $i$ choose one
$F^i_{k_i}$ equivalence class $X_i \sub Y_i$. We may assume that
either $f \rest \cup_i X_i$ is either constant or assumes distinct
values for different $i$'s, thus $\la \la F^i_j:j<k_i+1 ;X_i \ra: i
\in \omega \ra$ is the required extension.

\noi Otherwise $\limsup_i (m_i-k_i) = \infty$ and choose $X_i \sub
Y_i$ containing exactly one element from each $F^i_{k_i}$ equivalence
class. Then $\la \la F^i_j:k_i \leq j <m_i ;X_i \ra: i \in \omega \ra$
is now such that $f \rest X_i$ is one-one. It is now routine to
further extend the condition so that $f \rest \cup_i X_i$ is
one-one. This completes the proof. \qed

Thus restricted to some members of our ultrafilter, there will essentially
be only three kinds of functions in $\fo$; there is a
corresponding result for functions in $\fw$.

\begin{cor} \label{ww}
Given $h \in \fw$, and $\la \la F^i_j:j<m_i ;Y_i \ra: i \in \omega \ra
\in {\Bbb P}$, then there is   
\[ \la \la E^i_j:j<n_i ;X_i \ra: i \in \omega \ra 
   \leq  \la \la F^i_j:j<m_i ;Y_i \ra: i \in \omega \ra \]

\noi such that the order type of $h''\cup_i X_i$ is at most
$\omega^\omega$. 
\end{cor}

\noi To conclude the proof of Theorem \ref{simple}, let ${\Bbb G}$ be
a generic filter on ${\Bbb P}$, and $\CU$ the filter generated by

\[ \{\cup_i X_i :  \la \la E^i_j:j<n_i ;X_i \ra: i \in \omega \ra \in
{\Bbb G} \}.\]

\noi By Lemma \ref{c-closed}, every $X \sub \omega$ belongs to the
ground model, and by Lemma \ref{fw} (by considering characeristic
functions), $\CU$ contains a set $Y$ either included or disjoint from
$X$; thus $\CU$ is an ultrafilter. The nature of $\CU$ implies that it
cannot be better than an $\omega^\omega$-ultrafilter and Lemma \ref{ww}
shows that in fact it is an $\omega^\omega$-ultrafiler. Lemma \ref{fw}
also shows that $\CU$ has only one RK-predecessor, necessarily a Ramsey
ultrafilter. \qed

\vspace{.5in}

\noi {\bf Proof of Theorem \ref{complex}} We use the following partial order.

\begin{defn}

${\Bbb Q} = \{ \la X, \la E_\beta(X): \be \leq \al \ra \ra : X \in \so, \al
< \omega_1 \}$ where each $E_\be(X)$is an equivalence relation on $X$
with infinitely many classes and for each finite subset
$\{\be_1,\be_2,\dots ,\be_n \}$ of $\al$ (listed in increasing order)
the sequence $\la E_{\be_1}(X),$ $ E_{\be_2}(X),\dots, E_{\be_n}(X),
E_\al(X) \ra$ is \eifiner.  We further assume to simplify notation that
$E_0(X)$ is the trivial relation, equality.

\end{defn}

\noi We define the ordering as follows:

\[ \la X, \la E_\beta(X): \be \leq \al \ra \ra  \leq
             \la Y, \la E_\beta(Y): \be \leq \ga \ra \ra  \]

\noi if and only if $\ga \leq \al$ and for each $\be \leq \ga$, for all but
finitely many $E_\al(X)$ equivalence classes $C$, 
  $E_\be(X) \rest C = E_\be (Y) \rest C$.

\noi One should quickly verify that this indeed defines a transitive
ordering.

\begin{lem}\label{q-c-closed}
${\Bbb Q}$ is countably closed.
\end{lem}

\proof Given a decreasing sequence

\[ \la X_{n+1}, \la E_\beta(X_{n+1}): \be \leq \al_{n+1} \ra \ra  \leq
             \la X_n, \la E_\beta(X_n): \be \leq \al_n \ra \ra , \]

\noi for each $n \in \omega$ where we may as well assume that the
$\al_n$'s are strictly increasing, we let $\al = \sup_n \al_n$ and
construct

\[ \la X, \la E_\beta(X): \be \leq \al \ra \ra  \leq
             \la X_n, \la E_\beta(X_n): \be \leq \al_n \ra \ra , \]
	      
\noi for each $n$ as follows.  List $\al = \{\de_k : k \in \omega \}$
and proceed in $\omega$ steps to define the $E_\al$ equivalence
classes $\la E^i_\al(X): i \in \omega \ra$ on $X$ and thus $X$
itself.

\noi Having already defined the classes $E^j_\al(X)$ for $j<i$, choose
$n_i$ large enough so that $\{\de_j : j < i \} \sub \al_{n_i}$ and
$n_i > \{n_j: j < i \}$, and choose an $E_{\al_n}(X_n)$ equivalence
class $C$ for which the sequence $\{E_{\de_j}(X_n): j < i \}$, when
listed in increasing order of indices, is \ifiner on $C \sm \cup_{j<i}
E^j_\al(X)$, and such that $E_{\de_j}(X_n) \rest C = E_{\de_j}(X_j)
\rest C$ for all $j<n$.  Now simply let $E^i_\al(X) = C \sm \cup_{j<i}
E^j_\al(X)$. For $j \geq i$, we can define $E_{\delta_j}(X)$
arbitrarily on $E^i_\al(X)$.  \qed

\begin{lem}\label{q-shrink}

Given $f \in \fo$ and $\la Y, \la E_\beta(Y): \be \leq \al \ra \ra \in
{\Bbb Q}$, there is $ \la X, \la E_\beta(X): \be \leq \al \ra \ra \leq \la
Y, \la E_\beta(Y): \be \leq \al \ra \ra $ such that $f$ is either constant
on $X$ or else there is $\be \leq \al$ such that $f$ is constant on the
$E_\be(X)$ equivalence classes but assumes different values for
different classes. 

\end{lem}

\proof List $\al = \{ \de_k : k \in \omega \}$ and we may as well
assume that $\{ E_{\de_i}(Y) : i \leq k \} \cup \{E_\al(Y) \}$ is
\ifiner (listed in increasing order of indices) when restricted to the
$k^{th}$ class $E^k_\al(Y)$.  We may also assume that for each such
$k$ there is a $\be_k \in \{\de_i : i \leq k \}$ such that $f \rest
E^k_\al(Y)$ is constant on the $E_{\be_k}$ classes.  If $\be = \sup_k
\be_k$, we can further shrink $Y$ so that $\be_k=\be$ for all $k$.
When this process cannot yield a greater value for $\be$, then we can
require that $f$ assumes disitinct values for distinct $E_\be$
classes, this is the desired $X$.  \qed

\noi To conclude the proof of Theorem \ref{complex}, let ${\Bbb G}$ be
a generic filter on ${\Bbb Q}$, and $\CU$ the ultrafilter generated by

\[ \{ X \in \so : \la X, \la E_\beta(X): \be \leq \al \ra \ra \in {\Bbb G} \} \]

\begin{lem}\label{proper}

$\CU$ is a proper $\omega^\omega$-ultrafilter. 

\end{lem}

\proof By considering characteristic functions and using Lemma
\ref{q-shrink}, $\CU$ is an ultrafilter. Now let $\la X, \la E_\be(X): \be
\leq \al \ra \ra \in {\Bbb Q}$, list $\al = \{ \de_k : k \in \omega
\}$ and we assume again that $\{ E_{\de_i}(X) : i \leq k \} \cup
\{E_\al(X) \}$ is \ifiner (listed in increasing order of indices) when
restricted to the $k^{th}$ class $E^k_\al(X)$.

We first show that every function $h \in \fw$ can be restricted to a set
$X \in \CU$ so that its range has order type at most $\omega^\omega$.  For
this it suffices to shrink each $E_\al$ class so that actually the order
type of the range of $h$ restricted to the $E^k_\al(X)$ class is
at most $\omega^{k+1}$ and lies entirely after the range of $h$
restricted to the previous classes.  But then the order type of the
range of $h$ is at most $\omega^\omega$ as desired. 

We finally show that $\CU$ is a proper $\omega^\omega$-ultrafilter by
constructing an $h \in \fw$ whose range restricted to members of $\CU$
never drops below $\omega^\omega$.  With $\la X, \la E_b(X): \be \leq
\al \ra \ra $ as above, define $h$ as follows.  Let $\{\de_i : i \leq
k \} \cup \{\al \}$ be listed in increasing order as $\la \be^k_i : i
\leq k+1 \ra$ (so $\be^k_{k+1}=\al$).  We have by assumption that $\la
E_{\be^k_i}(X):i\leq k+1 \ra$ is \ifiner restricted to $E^k_\al(X)$.
Similarly to Lemma \ref{infchain}, define $h$ such that for each
$E_{\be^k_{i+1}}$ class, if the $E_{\be^k_i}$ subclasses are listed in
a sequence $E^\ell_{\be^k_i}$, then the range restricted to
$E^\ell_{\be^k_i}$ precedes the range restricted to $E^{\ell +
1}_{\be^k_i}$.  We may as well define $h$ to be constant on the
$E_{\be^k_0}$ classes.  Thus the order type of $h''X$ is $\oo$.

\noi Now if $ \la Y, \la E_\beta(Y): \be \leq \ga \ra \ra \leq \la X,
\la E_\beta(X): \be \leq \al \ra \ra $, choose for each $k$ an
$E_\ga(Y)$ class $C$ on which $\{E_{\be^k_i}(Y):i\leq k+1 \} \cup
\{E_\ga\}$ is \ifiner and $E_{\beta^k_i}(Y) \rest C = E_{\beta^k_i}(X)
\rest C $.  Then the range of $h$ restricted to this class has order
type at least $\omega^{k+1}$, and thus the order type of $h''Y$ is at
least $\oo$. \qed

Finally, by Lemma \ref{q-shrink}, every RK-predecessor of $\CU$ is
itself $\Bbb Q$-generic and therefore again a proper
$\omega^\omega$-ultrafilter by Lemma \ref{proper}.  This concludes the proof of
the theorem. \qed

\section{Conclusion}

It is  a natural step to consider next
$\omega^{\omega+\omega}$-ultrafilters and one interesting from
\cite{Bau} remains:

\noi {\bf Open Problem 3:} Does every
$\omega^{\omega+\omega}$-ultrafilter has an $\oo$ RK predecessor?

\noi The point is that for an $\omega^\alpha$-ultrafilter to have all
its RK predecessors also $\omega^\alpha$-ultrafilters, then $\alpha$
must also be indecomposable. Actually it is not hard to realize that
an $\omega^{\omega+\omega}$-ultrafilter must have a RK predecessor at
most an $\oo$-ultrafilter. The question is thus whether we can bypass
the value $\oo$.

\bibliographystyle{amsplain}

\end{document}